\documentclass[12pt]{amsart}

\usepackage{color}
\usepackage{amsfonts}
\usepackage{amscd}
\usepackage{amssymb}
\usepackage{amsthm}
\usepackage{amsmath}
\usepackage{comment}
\usepackage{epsfig}
\usepackage[all]{xy}

\usepackage[colorlinks=true, pdfstartview=FitV, linkcolor=blue, citecolor=blue, urlcolor=blue]{hyperref}
\usepackage{tikz,etex,bbm,xspace,hyperref}
\newcommand{\arxiv}[1]{\href{http://arxiv.org/abs/#1}{\tt arXiv:\nolinkurl{#1}}}

\newtheorem{Theorem}{Theorem}[section]

\newtheorem{Nested Lemma}[Theorem]{Nested Lemma}

\theoremstyle{definition}

\newtheorem{Definition}[Theorem]{Definition}

\theoremstyle{remark}
\newtheorem*{Remark}{Remark}

\newcommand{\Id}{\operatorname{Id}}

\makeatletter
\renewcommand{\@makefnmark}{\mbox{\textsuperscript{}}}
\makeatother

\newcommand{\h}{\hbar}

\renewcommand{\O}{\mathcal{O}}

\newcommand{\C}{\mathbb{C}}

\newcommand{\e}{\tilde{e}}

\newcommand{\nc}{\newcommand}
\nc{\gl}{\mathfrak{gl}}
\nc{\cosoc}{\operatorname{cosoc}}
\nc{\soc}{\operatorname{soc}}
\nc{\asl}{\widehat{\mathfrak{sl}}}
\nc{\g}{\mathfrak{g}}
\nc{\br}{\Bbb R}
\nc{\bz}{\Bbb Z}
\nc{\bn}{\Bbb N}
\nc{\omu}{\overline{\mu}}
\nc{\olambda}{\overline{\lambda}}
\nc{\oa}{\overline{a}}
\nc{\Irr}{\text{Irr }}
\nc{\oT}{\overline{T}}
\nc{\oR}{\overline{R}}
\nc{\oI}{\overline{I}}
\nc{\bfv}{{\bf v}}
\nc{\vareps}{\varepsilon}
\nc{\MV}{\mathcal{MV}}
\nc{\be}{\beta}
\nc{\End}{\text{End}}

\begin{document}

\title[Peter Weyl bases and preferred deformations]{Peter-Weyl bases, preferred deformations, and Schur-Weyl duality}

\author{Anthony Giaquinto}
\address{Department of Mathematics and Statistics, Loyola University, Chicago, IL}
\email{agiaqui@luc.edu}

\author{Alex Gilman}
\address{School of Physics and Astronomy, University of Minnesota, Minneapolis MN 55455}
\email{agilman@physics.umn.edu}

\author{Peter Tingley}
\address{Department of Mathematics and Statistics, Loyola University, Chicago, IL}
\email{ptingley@luc.edu}

 \thanks{}
 
%

\begin{abstract}
We discuss the deformed function algebra $\O_{\hbar}(G)$ of a simply connected reductive Lie group $G$ over ${\Bbb C}$ using a basis consisting of matrix elements of finite dimensional representations. This leads to a preferred deformation, meaning one where the structure constants of comultiplication are unchanged. The structure constants of multiplication are controlled by quantum $3j$ symbols. We then discuss connections earlier work on preferred deformations that involved Schur-Weyl duality.

\end{abstract}

\dedicatory{Dedicated to Kolya Reshetikhin on the
occasion of his $60$th birthday}

\date{\today}
\maketitle

%

\section{Introduction}

Let $G$ be a connected reductive Lie group over ${\Bbb C}$, and let $\g$ be its Lie algebra. Associated to this data are two Hopf algebras, the commutative function algebra $\O(G)$ and the cocommutative universal enveloping algebra $U(\g)$. During the 1980s, various non-commutative and non-commutative quantizations of these Hopf algebras were independently introduced. The first example, now known as $U_{\h}(\mathfrak{sl}_2)$, was discovered by Kulish and Reshetikhin in \cite{KR81} in relation to the quantum inverse scattering method. Later, Drinfeld and Jimbo independently introduced the well studied quantized universal enveloping algebra $U_\h(\g)$. On the dual side, several approaches to quantizations of $\O(G)$ have been studied. The first was the quantum matrix bialgebra $\O_{\hbar}(M_2)$ introduced by Faddeev and Takhtajan in \cite{FT}, constructed using the monodromy matrix for the quantum Lax operator of the Liouville model. This approach was fully developed in the landmark work \cite{FRT} of Faddeev, Reshetikhin and Takhtajan in which a quantum Yang-Baxter $R$-matrix is used to deform the defining relations of the classical series of coordinate algebras $\O(G)$.

We should also mention a few other early approaches. In \cite{W87a, W87b}, Woronowicz developed the theory of compact quantum groups in the $C^*$-algebra framework by introducing the quantization $SU_{\mu}(2)$ in which the parameter $\mu$ is a positive real number. Matrix coefficients of finite dimensional representations play a key role in this theory. Another approach due to Manin \cite {M} constructs quantum coendomorphism bialgebras as universal objects coacting on a pair of quantum linear spaces.

Since $U(\g)$ is rigid as an algebra and $\O (G)$ is rigid as a coalgebra, the fact that
$U_{\hbar}(\g)$ and $\O_{\hbar}(G)$ are formal deformations implies that their finite dimensional representations and corepresentations correspond exactly to those for $U(\g)$ and $\O(G)$ respectively. In particular, these categories are semi-simple/cosemisimple. With this in mind, a dual approach may be taken by first studying the monoidal categories of corepresentations of $\O_{\hbar}(G)$ and $\O(G)$. Once enough is known about these categories, one can follow the generalized Tannaka-Krein theory to reconstruct the Hopf algebras, which must necessarily be isomorphic as coalgebras, see \cite{JS}.

Focusing more sharply on our main point, it is a natural question to find a so-called ``preferred" presentations of $U_{\hbar}(\g)$ and $\O_{\hbar}(G)$, where the algebra structure is completely unchanged for $U_{\hbar}(\g)$ and the coalgebra structure is completely unchanged for $\O_{\hbar}(G)$. 
With the usual generators and relations descriptions this seems to be hard since, for the natural bases, all structures are varying.

The purpose of this note is to discuss a preferred presentation for $\O_{\hbar}(G)$. The starting point is to view $\O_{\hbar}(G)$ as the restricted dual Hopf algebra of $U_{\hbar}(\g)$. The preferred presentation is achieved from a Peter-Weyl basis of $\O_{\hbar}(G)$ -- a basis consisting of matrix elements of finite dimensional representations. The structure constants for the preferred presentation make use of quantum $3j$-symbols from physics,
which encode the decomposition of a tensor product of irreducible representations into
irreducibles. These coefficients have numerous applications and in the rank one case have been extensively studied, see \cite{KK89, KR88, V89}.

We finish by describing how this relates to Schur-Weyl duality in type A, and hence to some earlier work by Gerstenhaber, Giaquinto and Schack \cite{GGS:1992, Giaquinto:1992} on preferred deformations. In these papers, the quantum matrix bialgebra $\O _{\hbar}(M_n)$ is viewed as the invariant or ``quantum symmetric'' elements of the tensor algebra $T(M_n ^*)$ which are fixed by the action of a certain quantum symmetric group. This is a subgroup of the cactus group studied in e.g. \cite{KT}, and as discussed there is related to using Drinfeld's unitarized $R$-matrix from \cite{D90} in place of the usual $R$-matrix. If $V$ is the vector representation, then the image of this group in $\rm{End}(V^{\otimes n})$ generates the usual action of the Hecke algebra on this space.

In \cite{GGS:1992, Giaquinto:1992} the decomposition of tensor space $V^{\otimes n}$ into quantum symmetric elements is obtained with the aid of the Woronowicz quantization of $U(\mathfrak{sl}_n)$ which acts as skew derivations of $V^{\otimes n}$ associated to certain automorphisms. These automorphisms coincide with the exponentials $\exp (\h \, H_i)$ of the standard Cartan generators $H_i$ of $\mathfrak{sl}_n$. Thus the images of the Woronowicz quantization and $U_{\h}(\mathfrak{sl}_n)$ in $\rm{End}(V^{\otimes n})[[\h]]$ coincide and so the Schur-Weyl decompositions of $V^{\otimes n}$ are the same for either of these two quantizations. The use of the Woronowicz quantization was motivated by the fact that its finite dimensional representations correspond exactly to those of $U_\h(\mathfrak{sl}_n)$ or $U(\mathfrak{sl}_n)$, and one does not have to exclude the non-type-1 representation that appear for the rational form $U_{q}(\mathfrak{sl}_n)$. The disadvantage is that the Woronowicz quantization does not give a bialgebra structure (see \cite[p26]{GGS:1992}).

We do not carefully address the preferred presentation of $U_{\h}(\g)$ in this note. The dual Peter-Weyl basis does give a preferred presentation, but of a certain completion of $U_{\h}(\g)$. Finding a preferred presentation for the $U_{\h}(\g)$ itself seems more difficult. An algebra isomorphism from $U_{\h}(\mathfrak{sl}_2)$ to $U(\mathfrak{sl}_2)$ was given in \cite[Proposition 6.4.6]{CP}, giving the preferred presentation in that case. Only recently in \cite{AG17} was an explicit trivialization of $U_{\h}(\mathfrak{sl}_n)$ given by Appel and Gautam.  This isomorphism is induced by a map between the quantum loop algebra of $\mathfrak{sl}_n$ and a completion of the Yangian.

This note is organized as follows. In \S\ref{S:def} we discuss some background on deformation theory and the notion of preferred deformations. In \S\ref{S:PWdeform} we construct a preferred deformation using a Peter-Weyl basis. In \S\ref{S:SW1} we restrict to type $A$ and reformulate the construction using Schur-Weyl duality, then discuss how this relates to some older work.

\section{Deformations} \label{S:def}

\subsection{Formal deformations} Let $B$ be a bialgebra over $\C$. A $\C[[\h]]$-bialgebra $B_\h$ is a \emph{formal deformation} of $B$ if it is a topologically free $\C[[\h]]$-module together with an isomorphism $B_\h/hB_\h \simeq B$.
If $B_\h$ is a formal deformation of $B$, we can choose an identification of $B_\h$ with $B[[\h]]$ as a $\C[[\h]]$-module. The multiplication $\mu$ and comultiplication $\Delta$ of $B_\h$ then necessarily have the form
\begin{eqnarray*}
\mu(a,b)&=& \mu_0(a,b)+ h\mu_1(a,b)+h^2\mu_2(a,b)+\cdots \\
\Delta(a)&=& \Delta_0(a)+\h\Delta_1(a)+h^2\Delta_2(a)+\cdots
\end{eqnarray*}
where $\mu_i:A\otimes A\to A$ and $\Delta_i:A\to A\otimes A$ are linear maps and $\mu_0$ and $\Delta_0$ are the undeformed multiplication and comultiplication of $B$.

\subsection{Equivalence and preferred presentations} Deformations $B_\h$ and $B_\h'$ are equivalent if there is a ${\Bbb C}[[\h]]$-bialgebra isomorphism $\phi:B_\h\to B_\h'$ which reduces to the identity modulo $\h$.
It is known that \emph{every} deformation of the universal enveloping algebra $U(\g)$ is equivalent to one in which $\mu_\h=\mu$. That is, it is a trivial deformation of the algebra structure, and so the representation theory of $U_\h(\g)$ and $U(\g)$ is identical. Dually, every deformation of $\O(G)$ is equivalent to one in which $\Delta_\h=\Delta$, so the co-representation theory is unchanged. A \emph{preferred presentation} of a deformation of $U(\g)$ or $\O(G)$ is one with unchanged multiplication or comultiplication.

A natural question is to find preferred presentations of $U_\h(\g)$ and $\O_\h(G)$. To do so requires an identification of their underlying vector spaces with $U(\g)[[\h]]$ and $\O(G)[[\h]]$ as $\C[[\h]]$-modules. This can be accomplished, for example, by finding bases of $U_\h(\g)$ and $\O_\h(G)$ which reduce to bases of $U(\g)$ and $\O(G)$ modulo $\h$. However, most choices of bases do not provide preferred presentations: both multiplication and comultiplication depend on $\h$. This is true, in particular, for the various PBW-type bases in the literature.

We now arrive at an interesting juncture: Once $\O_\h(G)$ is shown to be a formal deformation, we know its irreducible representations are the same as those for $\O(G)$, just tensored with ${\Bbb C}[[\h]]$. We can then consider a Peter-Weyl type basis of $\O_\h(G)$, meaning a basis consisting of matrix elements of irreducible representations.
We shall see that this provides the sought after preferred presentation of $\O_\h(G)$.

\subsection{Standard deformation of $U(\g)$} \
Consider the standard Chevalley generators  $E_i,F_i,H_i$ of $U(\g)$. The deformation $U_\h(\g)$ is usually defined to be the algebra with these same generators, but deformed relations and a deformed coproduct.
The main structure we will need here is the coproduct, so we state that explicitly:
$$
\begin{aligned}
& \Delta E_i= E_i \otimes e^{-h H_i} + 1 \otimes E_i \\
& \Delta F_i= F_i \otimes 1 + e^{h H_i}  \otimes F_i \\
& \Delta H_i= H_i \otimes 1+ 1 \otimes H_i
\end{aligned}
$$
The rest of the structure can be found in many places, see e.g. \cite{CP}.
The relations for multiplication are also deformed. For instance, in $\mathfrak{sl}_2$ the undeformed relation $EF-FE=2H$ becomes
$$EF-FE=\frac{e^{hH}-e^{-hH}}{e^\h-e^{-\h}}.$$
So the deformation is certainly not preferred with respect to any PBW type basis.


%

\subsection{Standard deformation of $\O(M_n)$} \label{SS:OM-standard}

Here we give the deformed relations for $\O(M_n)$, as constructed by Faddeev-Reshetikhin-Takhtajan \cite{FRT}. We focus on the case $n=2$, and for simplicity of presentation define $q=e^\h$. Consider the coordinate functions
$$X=\begin{bmatrix}a&b\\c&d\end{bmatrix}=\begin{bmatrix}e_{11}^*&e_{12}^*\\e_{21}^*&e_{22}^*\end{bmatrix}$$
on the space of $2\times 2$ matrices. The FRT formalism starts with a solution $R$ to the quantum Yang-Baxter equation ($R_{12}R_{13}R_{23}=R_{23}R_{13}R_{12}$) and imposes the relations $RX_1X_2=X_2X_1R$ where $X_1=X\otimes I$ and $X_2=I\otimes X$. In coordinates,
$$
X_1X_2=\begin{bmatrix}a^2&ab&ba&b^2\\ac&ad&bc&bd\\ca&cb&da&db\\ c^2&cd&dc&d^2\end{bmatrix} , \quad
X_2X_1= \begin{bmatrix}a^2&ba&ab&b^2\\ca&da&cb&db\\ac&bc&ad&bd\\c^2&dc&cd&d^2 \end{bmatrix}.$$
For $M_2$,
$$R=\begin{bmatrix}q&0&0&0\\0&1&0&0\\0&q-q^{-1}&1&0\\0&0&0&q\end{bmatrix}$$
which produces the relations
\begin{equation} \label{eq:Rquant-rels}
\begin{aligned}
&ab=qba, \quad ac=qca, \quad bd = q db, \quad cd = q dc, \\
& \mbox{} \hspace{1cm} bc=cb, \quad ad-da=(q-q^{-1})bc.
\end{aligned}
\end{equation}
The coproduct is defined on generators by
$$\Delta \left(\begin{bmatrix}a&b\\c&d\end{bmatrix}\right)= \begin{bmatrix}a&b\\c&d\end{bmatrix} \otimes \begin{bmatrix}a&b\\c&d\end{bmatrix}
=
\begin{bmatrix}a \otimes a+ b \otimes c &a \otimes b+ b \otimes d \\c \otimes a + d \otimes c &c \otimes b+ d \otimes d\end{bmatrix}
$$
and is extended multiplicatively to monomials of higher degree.
For example
$$
\Delta (a^2) =
a^2\otimes a^2+ ab \otimes ac + ba \otimes ca+b^2 \otimes c^2 = a^2\otimes a^2  +(1+q^{-2})ab\otimes bc +b^2\otimes c^2.
$$
This is dependent on $q$, so the deformation is not preferred, at least when using the PBW basis $\{a^i\,b^j\,c^k\,d^l\quad |\quad i,j,k,l \in \mathbb Z_{\geq 0}\}$ to identify $\O_\h(M_n)$ with $\O(M_n) [[\h]]$.


\section{Peter-Weyl bases and preferred deformations} \label{S:PWdeform}

Another approach to deforming $\O(G)$ is by duality: one simply defines $\O_{\h}(G)$ as the restricted dual of $U_{\h}(\g)$. In this setting the Peter-Weyl basis arises naturally. It is with this basis that we get a preferred presentation of $\O_\h(G)$.

\subsection{The Peter Weyl basis}
Since $\O_\h(G)$ is cosemisimple, the restricted dual definition implies
\begin{equation*}
\O_\h(G) \simeq \oplus_\lambda \End (V_\lambda)^*,
\end{equation*}
where the $\lambda$ runs over the dominant integral weights of $G$, the $V_\lambda$ are the corresponding representations of $U_\h(\mathfrak{g})$, and the isomorphism is as coalgebras over $k[[\h]]$. See e.g. \cite[Chapter 11]{KS97}. 
This becomes an isomorphism of Hopf algebras if one defines multiplication on $\oplus_\lambda \End (V_\lambda)$ as the dual of the coproduct for $U_{\h}(\g)$.

For any $\lambda$, $\End (V_\lambda)$ is naturally identified with $V_\lambda \otimes V_\lambda^*$. Taking duals, we identify  $\End (V_\lambda)^*$  with $V_\lambda^* \otimes V_\lambda$.
This gives a natural way to choose a basis for $\O_\h(G)$:
pick dual bases $B_\lambda$ and $B^*_\lambda$ for each pair $V_\lambda$,$V_\lambda^*$. Then
\begin{equation*} 
\bigsqcup_\lambda \{ Y^* \otimes X : X,Y \in B_\lambda\}.
\end{equation*}
is a basis for $\O_\h(G)$, which we call a {\bf Peter-Weyl basis}.
The pairing of $U(\g)$ with $\O(G)$ is given by,
for $Y^* \otimes X \in \O(G)$ and $u \in U(\g)$,
$$\langle Y^* \otimes X, u \rangle = Y^* u(X).$$

\begin{Remark}
One often reverses order of factors when taking duals of tensor products. We have not done so, in part to match conventions in \cite{GGS:1992}.
\end{Remark}

\subsection{Comultiplication}  \label{ss:comult} We are identifying $\prod_\lambda \End V_\lambda$ with a completion of $U_\h(\g)$, and $\O_\h(G)$ with the dual Hopf algebra to this. Thus, in both $\O(G)$ and $\O_\h(G)$,
co-multiplication is the dual of multiplication in $\prod_\lambda \End V_\lambda$. In coordinates, for $X, Y \in B_\lambda$,
\begin{equation}
\label{eq:comult}
\Delta(Y^* \otimes X )= \sum_{Z \in B_\lambda}  (Y^* \otimes Z) \otimes (Z^* \otimes X).
\end{equation}

\subsection{Multiplication (abstract)} \label{ss:mult}
Multiplication is the dual to comultiplication in $U_\h(\g)$.
 In coordinates this means,  for  $X_1, Y_1 \in B_\lambda$, $X_2,Y_2 \in B_\mu$, and any $u \in U(\g)$,
\begin{equation}
\label{eq:mult1}
((Y_1^*\otimes X_1 )(Y_2^* \otimes X_2))(u)=  (Y_1^* \otimes Y^*_2) \Delta u (X_1 \otimes X_2).
\end{equation}
To be explicit we need to express this in terms of the Peter Weyl basis. The resulting structure constants are closely related to the famous 3j symbols from physics.

\subsection{$3j$ symbols}
These are often studied just for $\text{SL}(2)$, but we need the following more general notion.
For each triple $\lambda, \mu, \nu$, choose a basis $\{ \phi_1, \cdots, \phi_{c_{\lambda, \mu}^\nu} \}$ for the space of embeddings of $V_\nu \hookrightarrow V_\lambda \otimes V_\mu$. For $X_1 \in B_\lambda, X_2  \in B_\mu, X_3 \in B_\nu$, write
\begin{equation*}
X_1 \otimes X_2= \sum_\nu \sum_{1 \leq k \leq c_{\lambda,\mu}^\nu} \sum_{X_3 \in B_\nu } \left(
\begin{array}{ccc}
\lambda & \mu & \nu \\
X_1 & X_2 & X_3
\end{array}
 \right)_k
 \phi_k(X_3).
\end{equation*}
The constants
$\displaystyle  \left(
\begin{array}{ccc}
\lambda & \mu & \nu \\
X_1 & X_2 & X_3
\end{array}
 \right)_k$
are called the $3j$ symbols.

Taking duals gives a basis
 $\{ \phi^*_1, \cdots, \phi^*_{c_{\lambda, \mu}^\nu} \}$  of the space of surjections $V^*_\lambda \otimes V^*_\mu \rightarrow V_\nu^*$. We then get dual $3j $ symbols defined by
$$
\phi_k^*(Y_1^* \otimes Y_2^*)= \sum_{Y_3^* \in B_\nu^*}
 \overline{\left(
\begin{array}{ccc}
\lambda & \mu & \nu \\
Y_1^*& Y_2^* & Y_3^*
\end{array}
 \right)}_k Y_3^*.
 $$

This can be done just as easily for representations of $U(\g)$ or $U_{\h}(\g)$.

\begin{Remark}
 In the $\text{SL}(2)$ case, there is a unique (up to signs) orthonormal weight basis, so a chosen Peter-Weyl basis. The spaces of embeddings $V_\nu \hookrightarrow V_\lambda \otimes V_\mu$ are 1 dimensional, and the inner product can be used to  normalize the embedding, fixing the $3j$ symbols. As mentioned earlier, these have been calculated extensively. Using orthonormal bases also implies that the the $3j$ symbols and dual $3j$ symbols coincide exactly.
\end{Remark}


 \subsection{Structure constants for multiplication} \label{ss:sc}
 It is now immediate from definitions that the structure constants for multiplication in the Peter-Weyl basis are given by, for  $X_1, Y_1 \in B_\lambda$ and $X_2,Y_2 \in B_\mu$,
 \begin{equation}
 \label{eq:structure-constants}
 (Y_1^* \otimes X_1) (Y_2^* \otimes X_2)   \hspace{-0.1cm} =
 \hspace{-0.5cm}
 \sum_{ \substack{\nu \in P_+ \\ {{X_3, Y_3 \in B_\nu}}}}
  \hspace{-0.3cm}
 \left[
 \sum_{1 \leq k \leq c_{\lambda, \mu}^\nu}
 \overline{
  \left(
\begin{array}{ccc}
\lambda & \mu & \nu \\
Y_1^* & Y_2^* & Y_3^*
\end{array}
 \right)}_k
  \left(
\begin{array}{ccc}
\lambda & \mu & \nu \\
X_1 & X_2 & X_3
\end{array}
 \right)_k
 \right]
  \hspace{-0.1cm}
  Y_3^* \otimes X_3.
 \end{equation}
Multiplication as defined in \eqref{eq:mult1} does not depend on the basis $\{\phi_1, \ldots, \phi_k\}$, so
 $$
 \sum_{1 \leq k \leq c_{\lambda, \mu}^\nu}
 \overline{
  \left(
\begin{array}{ccc}
\lambda & \mu & \nu \\
Y_1^* & Y_2^* & Y_3^*
\end{array}
 \right)}_k
  \left(
\begin{array}{ccc}
\lambda & \mu & \nu \\
X_1 & X_2 & X_3
\end{array}
 \right)_k
 $$
 must be independent of the choice of basis $\{\phi_1, \ldots, \phi_k\}$ as well.


 \subsection{Preferred presentation of $\O_\h(G)$}
 The set of irreducible representations $V_\lambda$ of $U(\g)$ and $U_\h(\g)$ correspond exactly, and
the comultiplication from \S\ref{ss:comult} does not reference $U(\g)$ at all, so is unchanged under deformation.
 The multiplication from \S\ref{ss:mult} does change when we move to $U_\h(\g)$, since it's definition uses the coproduct of $U(\g)$, which is deformed in $U_\h(\g)$. But the presentation in \S\ref{ss:sc} is still valid. The only difference is that the spaces of embeddings $V_\nu \hookrightarrow V_\lambda \otimes V_\mu$ change.

 In order to see $\O_\h(\g)$ as a preferred deformation of $\O(\g)$, one must simply choose
 \begin{itemize}
 \item A basis for each $U_\h(\g)$ module $V_\lambda$ which specializes to a basis at $h=0$,

 \item A basis for each space of $U_q(\g)$-homomorphisms  $V_\nu \hookrightarrow V_\lambda \otimes V_\mu$ which also specializes to a basis at $\h=0$. This leads to a definition of quantum 3j symbols.
 \end{itemize}
 Then the construction above gives a deformation where the structure constants of comultiplication are manifestly identical, and the structure constants for multiplication are given by \eqref{eq:structure-constants}, but with the 3j symbols replaces by their deformed counterparts.


\begin{Remark}
We have relied on the fact that we already have a (non-preferred) deformation of $U(\g)$ to construct our preferred presentation of $\O_\h(G)$.
One might try to use this approach to construct a deformation from scratch, by simply deforming the spaces of embeddings  $V_\nu \hookrightarrow V_\lambda \otimes V_\mu$. However, this deformation is not arbitrary: one needs to ensure that $\O(G)$ remains a Hopf algebra. Directly ensuring this seems difficult.
\end{Remark}


\subsection{Preferred deformation of $U(\g)$}
This method also gives a preferred presentation of a completion of $U_\h(\g)$ by working with the topological basis $\{ f_{b,c}^\lambda = b \otimes c^* \in \End V_\lambda\}$. The operations are the duals those of $\O_\h(G)$.
However, $U_\h(\g)$ is a proper subalgebra of $\prod_\lambda \End(V_\lambda)$, and the preferred presentation does not restrict in any nice way. It is also unclear how to relate the Peter-Weyl type bases with the Chevalley generators. So this approach is not really satisfactory.

\subsection{Non-simply-connected groups and matrix algebras}

The condition that $G$ be simply connected is not really needed. A non-simply connected reductive Lie group $G'$ is always the quotient of a corresponding simply-connected one, and the category of finite dimensional representations of $G'$ is a sub-tensor-category of the category of finite dimension representations of $G$. The irreducible representations of that category are parameterized by $\lambda$ in the positive part of some sub-lattice $P'$ of the weight lattice of $G$. The whole story then goes through by realizing $\O(G')$ as $\oplus \End(V_\lambda)^*,$ where now one restricts to $\lambda \in P'_+$.

In type $A$ one can also consider $\mathcal{O}(\text{GL}_k)$, and again the story goes through without significant changes, only now there are more representations than for  $\mathcal{O}(\text{SL}_k)$, since any irreducible representation can be tensored with any integer power of the determinant representation.
In \S\ref{S:SW1} we will actually work with $\mathcal{O}(M_k)$, the function algebra on the algebra of all $k \times k$ matrices. This is isomorphic to $\oplus_\lambda \End V_\lambda^*$, where now the $\lambda$ index the polynomial representations of $\text{GL}_k$. These $\lambda$'s are naturally indexed by partitions with at most $k$ parts.


\section{Relation to Schur-Weyl Duality} \label{S:SW1}

 For the case of $\rm{GL}_k$, or $\rm{SL_k}$, \cite{Giaquinto:1992, GGS:1992} studied another approach to finding a preferred deformation. Their approach most naturally realizes $\O(M_k)$, the function algebra on the algebraic monoid of $k\times k$ matrices. We now discuss how their results naturally arise in our framework.

 \subsection{General categorical discussion}
Identifying $\mathcal{O}(G)$ with $\oplus_\lambda \End(V_\lambda)^*$ as we have is not necessarily the most natural thing to do, since it requires choosing a representation in each isomorphism class of simples. In fact, any $f \in \End(V)^*$, for any representation $V$, gives a function of $G$. However different elements of $\End(V)^*$ can give identical functions on $G$, so $\mathcal{O}(G)$ should be identified with a quotient of $\oplus_V \End(V)^*$. This is a badly infinite sum, but ignoring that for now, multiplication is simple: given two elements of $\mathcal{O}(G)$, $f \in \End(V)^*, g \in \End(W)^*$, $$fg = f \otimes g \in \End(V \otimes W)^* \simeq \End(V)^* \otimes \End (W)^*.$$
We work with $\oplus_\lambda (\End V_\lambda)^*$ essentially because every element of $\mathcal O(G)$ appears exactly once. Equivalently, using the restricted dual definition, every linear functional on $U(\g)$ is represented only once.


 \subsection{Using $V^{\otimes n}$, undeformed}  \label{ss:undeformedVnrealization}
 Now we will restrict to considering $\mathcal{O}(M_k)$. Then there is another natural space which encodes every function exactly once:
 $$\oplus_n ( (\End V^{\otimes n})^{S_n})^*,$$
 where $V$ is the vector representation.
 To see this, recall that, by Schur-Weyl duality,
 $$  V^{\otimes n} \simeq  \oplus_\lambda V_\lambda \boxtimes W_\lambda,$$
  where $\lambda$ ranges over all partitions of $n$ with at most $k$ rows, and $V_\lambda, W_\lambda$ are the irreducible representations of $M_k$ and $S_n$ respectively.
 Then, by Schur's lemma,
\begin{equation*} 
\oplus_n (\End V^{\otimes n})^{S_n} \simeq \oplus_\lambda \End V_\lambda.
\end{equation*}

 We will need to understand this identification explicitly. Fix $\lambda,$ and choose any $w,w^* \in W_\lambda, W_\lambda^*$ with $w^*(w)=1$. Then, for any
 $c^* \otimes b \in \End(V_\lambda)^*$, the element $(c^* \otimes w^*_m) \otimes (b  \otimes m_w) \in \End (V^{\otimes n})^*$ gives the same function on $U_q(\mathfrak{gl}_n)$ as $c^* \otimes b$, but it is not $S_n$ invariant. To fix that,
consider the Young symmetrizer $\displaystyle P= \frac{1}{n!} \sum_{\sigma \in S_n} \sigma$. Then
 $$ P ( (c^* \otimes m_w^*) \otimes (b \otimes w_m))$$
 is clearly $S_n$ invariant, and gives the same function on $U_q(\mathfrak{gl}_n)$. Explicitly,
 $$P ((b^* \otimes w_m^*) \otimes (c \otimes m_w))= \frac{1}{n!} \sum_{\sigma \in S_n} (b^* \otimes g w_m^*) \otimes (c \otimes g m_w).$$
Since there is only one $S_n$ invariant element corresponding to a given function, this is independent of the choice of $w, w^*$. 

Multiplication is then given by, for  $f \in  (\End V^{\otimes n})^{S_n}$ and $g \in  \End (V^{\otimes m})^{S_m}$,
 $$ f g=P_{n+m}(f \otimes g).$$
Comultiplication would normally be given by, for any dual bases $C$ and $C^*$ of $V^{\otimes n}$,
 \begin{equation*}
\Delta(Y^* \otimes X)= \sum_{Z \in C}  (Y^* \otimes Z) \otimes ( Z^* \otimes X).
\end{equation*}
However, while this is in $\End (V^{\otimes n})\otimes \End (V^{\otimes n})$, and corresponds to the correct element of $\mathcal{O}(M_k) \otimes \mathcal{O}(M_k) $, it is not in $(\End V^{\otimes n})^{S_n} \otimes (\End V^{\otimes n})^{S_n}$ . To fix this, apply the symmetrizer $P$ to each of the two factors to get something in the right space which gives the same function on $U_{\h}(\text{gl}_n)$. The correct definition becomes
 \begin{equation*}
\Delta(Y^* \otimes X)= ( P \otimes P )\sum_{Z \in C}  (Y^* \otimes Z) \otimes ( Z^* \otimes X).
\end{equation*}

 \subsection{Using $V^{\otimes n}$, deformed}
By quantum Schur-Weyl duality $V^{\otimes n} \simeq \oplus_\lambda V_\lambda \boxtimes W_\lambda$, where $V$ is now the vector representation of $U_{\h}(\mathfrak{gl}_n)$, the $V_\lambda$ are the polynomial representation of $U_\h(\mathfrak{gl}_n)$, and the $W_\lambda$ are the irreducible representations of the Hecke algebra $H_n$ corresponding to partitions with at most $k$ rows. Then by Schur's lemma,
$$\oplus_n ((\End(V^{\otimes n})^{H_n})^* \simeq \oplus_\lambda \End(V_\lambda)^* = \mathcal{O}(M_n),$$
where the $H_n$ superscript means $H_n$ equivariant functions. The space on the left can be written as
$\oplus ((V^*)^{\otimes n} \otimes  V^{\otimes n})^{H_n},$
where the $H_n$ still means equivariant.

 The key to understanding the operations in \S\ref{ss:undeformedVnrealization} was to understand the projection
$$(V^*)^n \otimes V^n \rightarrow ((V^*)^n \otimes V^n)^{S_n},$$
where $S_n$ acts simultaneously on the two factors. This was defined as the Young symmetrizer $P$ acting simultaneously on both factors, but to quantize we need a different characterization,
 since it is not clear how to have $H_n$ act on $(V^*)^n \otimes V^n$.

The crucial thing in the previous section is that $P$ acts on $\oplus_n (\End V^n)^*$ as a projection so that, for any  $\phi \in (V^*)^n \otimes V^n $, $\phi$ and $P(\phi)$ define the same function on $U(\mathfrak{gl}_k)$. In this form, there is no problem giving the deformed definition.

 \begin{Definition}
$\pi:   \oplus_n (V^*)^{\otimes n} \otimes  V^{\otimes n}\rightarrow \oplus_n ((V^*)^{\otimes n} \otimes  V^{\otimes n})^{H_n} $ is the unique projection such that, for any $\phi \in  (V^*)^{\otimes n} \otimes  V^{\otimes n}$ , $\phi$ and $\pi(\phi)$ define the same function on $U_\h(\g)$.
\end{Definition}

This induces a Hopf algebra structure on $\oplus ((V^*)^{\otimes n} \otimes  V^{\otimes n})^{H_n}$ because the subset of $\oplus ((V^*)^{\otimes n} \otimes  V^{\otimes n})$ consisting of elements that define the zero function on $U_\h(\g)$ is a Hopf ideal.
Multiplication and comultiplication on $\oplus ((V^*)^{\otimes n} \otimes  V^{\otimes n})^{H_n}$ are given by:
\begin{equation}
\label{eq:comult-Vn-q}
\Delta(Y^* \otimes X)= ( \pi \otimes \pi  )\sum_{Z \in C}  (Y^* \otimes Z) \otimes ( Z^* \otimes X),
\end{equation}
\begin{equation} \label{eq:mult-Vn-q}
fg = \pi (f \otimes g).
\end{equation}

\begin{Remark}

It would be nice to have a more explicit formula for $\pi$. In the case $n=2$ such a formula is known. As we shall see in secton \ref{previous}, the $H_2$-equivariant endomorphisms of $V^*\otimes V^*$ are determined by an involution $Q$. It follows that $\pi=\frac{1+Q}{2}$, see \cite{Giaquinto:1992, GGS:1992}.



In general one might try to replace $P$ with the $q$-symmetrizer from \cite{Gyoja:1986}.
This does give a natural analogue of $P$ acting on $V^{\otimes n}$, but we would need it to act on $(V^*)^{\otimes n} \otimes V^{\otimes n}$.
If the $T_i$ are the generators of the Hecke algebra, the appropriate action on $(V^*)^n$ should replace $T_i$ with $T_i^{-1}$, and these satisfy a different set of Hecke-algebra relations. So the Hecke algebra does not even naturally act on $(V^*)^{\otimes n} \otimes V^{\otimes n}$. 
In fact no symmetrizer that acts by simultaneous permutations in the $(V^*)^{\otimes n}$ and $V^{\otimes n}$ can work, as then the relation
$ad-da=(q-q^{-1})bc$  in
 \eqref{eq:Rquant-rels} would not be possible.
\end{Remark}

\subsection{Preferred presentation} \label{ss:PSWD}
We can now construct a preferred presentation of $\O_\h(M_n)$. 

\begin{itemize}

\item Fix the Schur-Weyl duality isomorphism  $V^{\otimes n} \simeq \oplus_\lambda V_\lambda \boxtimes W_\lambda$. Of course one then gets a corresponding isomorphism
 $(V^{\otimes n})^* \simeq \oplus_\lambda V_\lambda^* \boxtimes W_\lambda^*$.

\item Fix bases $B_\lambda$ for each $V_\lambda$, and $P_\lambda$ for each $W_\lambda$, and their dual bases
 $B_\lambda^*$ for each $V_\lambda^*$, and $P_\lambda^*$ for each $W_\lambda^*$, in such a way that all specialize at $h=0$.

\item Then
$\displaystyle \quad  \bigcup_\lambda  \left \{ X^\lambda_{b,c^*}:=\frac{1}{\dim W_\lambda} \sum_{p \in P_\lambda}   (c^* \boxtimes p^*) \otimes (b \boxtimes p)  : b \in B_\lambda, c^* \in B_\lambda^* \right\}$

\noindent is a basis for $\mathcal{O}_\h(M_n)$. As a function on $U_\h(\mathfrak{g})$, the element $ X^\lambda_{b,c^*}$ agrees with $c^* \otimes b \in (\End V_\lambda)^*$.

\end{itemize}
Since $ X^\lambda_{b,c^*}$ agrees with $c^* \otimes b \in (\End V_\lambda)^*$, the structure constants of multiplication and comultiplication in this basis must agree with \eqref{eq:comult} and \eqref{eq:structure-constants}. It is an interesting exercise to directly obtain these formulae from the new definitions of comultiplication \eqref{eq:comult-Vn-q} and multiplication \eqref{eq:mult-Vn-q}.

\subsection{Comparing with previous work}\label{previous}

We now compare the current approach with the ``method of quantum symmetry" from \cite{Giaquinto:1992, GGS:1992}. The starting point there is to view $\O(M_k)$ as the symmetric algebra  $SX=\bigoplus _{n\geq 0} (X^{\otimes n})^{S_n}$, where $X=V^*\otimes V$. To quantize, $S_n$ is replaced by a ``quantum symmetric group" $qS_n$ with generators $\tau_1,\ldots ,\tau _{n-1} $ and relations $\tau_i^2=\Id$ and $\tau_i\tau_j=\tau_j\tau_i$ if $|i-j|>1$.
Note that if the braid relations $\tau_i \tau_{i+1}\tau_i=\tau_{i+1}\tau_i\tau_{i+1}$ are added then we have the Artin presentation of $S_n$. As mentioned in the introduction, $qS_n$ is a subgroup of the cactus group.

 To describe the $qS_{n}$-action on $X^{\otimes n}$ we first deform the flip operator $\sigma: V^*\otimes V^*\rightarrow V^*\otimes V^*$ where $\sigma(\alpha \otimes \beta)=\beta \otimes \alpha$. Let $r=\sum_{i<j}e_{ij}\wedge e_{ji} = \frac{1}{2}\, (e_{ij}\otimes e_{ji}-e_{ji}\otimes e_{ij})$. This is the standard unitary solution to the modified classical Yang-Baxter equation associated to $\mathcal O_{\h}(M_k)$. Define an involution of $V^*\otimes V^*$ by $Q=(\exp(-\h r))\sigma (\exp{\h r})$. With this there is an action of $qS_n$ on $(V^*)^{\otimes n}$ where $\tau_i$ acts as $Q$ in factors $i$ and $i+1$ and the identity elsewhere. Taking duals there is a corresponding action on  $(V^*)^{\otimes d}$ and hence $qS_n$ acts diagonally on $(V^*)^{\otimes n} \otimes V^{\otimes n}=X^{\otimes n}$.

 One of the main results of \cite{Giaquinto:1992, GGS:1992} is that the set of invariant elements of the tensor algebra $TX$ is a bialgebra which is isomorphic to $\mathcal O_{\h}(M_k)$. Moreover, the comultiplication in $\bigoplus_{n\geq 0} (X^{\otimes n})^{qS_n}$ is independent of $\h$ and coincides with the usual comultiplication in $\O(M_k)=\bigoplus_{n\geq 0} (X^{\otimes n})^{ S_n}$. Thus this construction yields the desired preferred presentation of $\mathcal O_{\h}(M_k)$.

 This essentially coincides with our construction.
 Using the notation of \cite[\S10]{GGS:1992}, $(k\langle M(n)\rangle^*, \otimes)$ is naturally the tensor algebra of $V^* \otimes V$, which we identify with 
 $\oplus_n (\End V^{\otimes n})^*$, and think of as functions on $U_\h(\g)$.
 The space $\text{sk}_q \langle M(n)^* \rangle$ is generated by the images of the operators $\frac{1}{2}(Id -\tau_i)$ acting on $(V^* \otimes V)^{\otimes n}$, and these images are easily seen to define the zero function on $U_\h(\g)$. 
 So, the quotient in the top line of the diagram in \cite[ Theorem 10.8]{GGS:1992} is by a set of elements which are all the zero function on $U_\h(\g)$, and by comparing dimensions it agrees with our $\pi$.
Thus the comultiplication given in \cite{Giaquinto:1992, GGS:1992} coincides exactly with \eqref{eq:comult} and \eqref{eq:comult-Vn-q}, and the multiplication is described using the projection formula \eqref{eq:mult-Vn-q}.

The expression for the multiplicative structure constants in terms of $3j$ symbols is largely new to this paper, although the multiplication formulas for quantum linear spaces given in \cite{Giaquinto:1992, GGS:1992} can easily be expressed in the $3j$ symbol notation, and this in turn gives some of the structure constants for $\mathcal{O}(M_n)$. So this idea really dates to those papers as well.

%

\subsection{Deriving the $R$-matrix relations in $\mathcal{O}_q(M_2)$}

We now derive the last two relations in the FRT construction of $\mathcal{O}_q(M_2)$  (see  \S\ref{SS:OM-standard}) in our language (the others are simpler).
One could also see that the constructions agree by directly showing that the $\frac{1}{2}(\Id - \tau)$ action on $X \otimes X$ gives the FRT relations.

The variables $a,b,c,d$ in our language are
$$
a= e_1 \otimes e_1^*, \quad b= e_1 \otimes e_2^*, \quad c= e_2 \otimes e_1^*, \quad d= e_2 \otimes e_2^*. \quad
$$
As a representation of $U_\h(\mathfrak{gl}_2)$, $V \otimes V \simeq W \oplus T$, where $W$ is a three dimensional representation and $T$ is one dimensional. These have basis
$$ W: \{ e_1 \otimes e_1, e_2 \otimes e_1+q e_1 \otimes e_2, e_2 \otimes e_2 \}, \qquad T: \{ e_2 \otimes e_1-q^{-1} e_1 \otimes e_2 \}.$$
Let $\quad s= e_2 \otimes e_1+q e_1 \otimes e_2, \quad t=  e_2 \otimes e_1-q^{-1} e_1 \otimes e_2.$
Then $\{ s,t\}$ spans the 0 weight space of $V \otimes V$. Let $\{ s^*, t^*\}$ be the dual basis of this weight space.  Then
\begin{equation*}
\begin{aligned}
& e_1 \otimes e_2 = \frac{s-t}{q + q^{-1}}, \quad e_2 \otimes e_1 = \frac{q^{-1}s+q t}{q+q^{-1} }
\\
&e^*_1 \otimes e^*_2 =q s^*-q^{-1}t^*, \quad e_2^* \otimes e_1^* = s^*+t^*.
\end{aligned}
\end{equation*}
The Hecke algebra is the algebra of operators commuting with the action of $U_\h(\mathfrak{gl}_2)$, so it is spanned by the projections onto $W$ and $T$. Thus $s^* \otimes s$ and $t^*\otimes t$ are both $H_2$ equivariant. Both $s^* \otimes t$ and $t^* \otimes s$ are zero as functions on $U_\h(\mathfrak{gl}_2)$ by Schur's lemma, so these are both killed by $\pi$. Thus

\vspace{-0.25cm}

\begin{equation*}
\begin{aligned}
ad & =\pi ((e_1^* \otimes e_2^* ) \otimes (e_1 \otimes e_2) )  \\
&= \pi \left((q s^*-q^{-1}t^*) \otimes  \frac{s-t}{q + q^{-1}} \right) \\
&= \frac{q}{{q + q^{-1}}} s^* \otimes s + \frac{q^{-1}}{{q + q^{-1}}} t^* \otimes t
\end{aligned}
\end{equation*}
\begin{equation*}
\begin{aligned}
da & =\pi ((e_2^* \otimes e_1^* ) \otimes (e_2 \otimes e_1) )  \\
&= \pi \left((s^*+t^*) \otimes  \frac{q^{-1}s+q t}{q+q^{-1}}  \right)\\
&= \frac{q^{-1}}{q+q^{-1}} s^* \otimes s +  \frac{q}{q+q^{-1}}t^* \otimes t
\end{aligned}
\end{equation*}
\begin{equation*}
\begin{aligned}
bc & =\pi ((e_2^* \otimes e_1^* ) \otimes (e_1 \otimes e_2) )  \\
&= \pi \left((s^*+t^*) \otimes \frac{s-t}{q + q^{-1}} \right)\\
&= \frac{1}{q+q^{-1}} s^* \otimes s -  \frac{1}{q+q^{-1}}t^* \otimes t
\end{aligned}
\end{equation*}
\begin{equation*}
\begin{aligned}
cb & =\pi ((e_1^* \otimes e_2^* ) \otimes (e_2 \otimes e_1) )  \\
&= \pi \left( (q s^*-q^{-1}t^*)  \otimes  \frac{q^{-1}s+q t}{q+q^{-1}}  \right)\\
&= \frac{1}{q+q^{-1}} s^* \otimes s -  \frac{1}{q+q^{-1}}t^* \otimes t.
\end{aligned}
\end{equation*}
Now the relation $bc=cb$ is obvious, and $ad-da=(q-q^{-1})bc$ is a simple calculation.

\end{document}